\documentclass{svmult}
\usepackage{mathptmx}
\usepackage{helvet}
\usepackage{courier}
\usepackage{graphicx}
\usepackage{makeidx}
\usepackage{multicol}
\usepackage{footmisc}
\usepackage{epsfig}
  \usepackage{t1enc}
      \usepackage[latin1]{inputenc}
      \usepackage[english]{babel}
          \usepackage{latexsym}
  \usepackage{amssymb}
\usepackage{amsfonts}
\def\theequation{\arabic{section}.\arabic{equation}}

\newtheorem{thm}{Theorem}[section]

\newcommand{\AS}{{\rm a.s.}}

\begin{document}

\title*{ON BAHADUR-KIEFER TYPE PROCESSES FOR SUMS AND RENEWALS IN  
DEPENDENT CASES}
\author{Endre Cs\'aki\and Mikl\'os Cs\"org\H{o}}
\institute{Endre Cs\'aki \at Alfr\'ed R\'enyi Institute of Mathematics,
Hungarian Academy of Sciences,\\ Budapest, P.O.B. 127, H-1364, Hungary,\\
\email{csaki.endre@renyi.mta.hu}
\and
Mikl\'os Cs\"org\H{o} \at
School of Mathematics and Statistics, Carleton University,
Ottawa, Ontario, Canada K1S 5B6\\
\email{mcsorgo@math.carleton.ca}}

\renewcommand{\theequation}{\thesection.\arabic{equation}}
\newcommand{\dto}{\stackrel{{d}}{\rightarrow}}
\newtheorem{propo}{Proposition}[section]
\newtheorem{cjt}{Conjecture}[section]

\titlerunning{Bahadur-Kiefer processes}
\maketitle

\abstract
{We study the asymptotic behaviour of Bahadur-Kiefer processes that are 
generated by summing partial sums of (weakly or strongly dependent) random 
variables and their renewals. Known results for i.i.d. case will be extended 
to dependent cases.} 

\vspace{.1cm}



\section{Introduction}\label{sect1}
 \renewcommand{\thesection}{\arabic{section}} \setcounter{equation}{0}
 \setcounter{thm}{0} \setcounter{lemma}{0}

In this work we intend to deal with Bahadur-Kiefer type processes that
are based on partial sums and their renewals of weakly, as well 
as strongly, dependent sequences of random variables. In order to initiate
our approach, let $\{Y_0,Y_1,Y_2,\ldots\}$ be random variables which have the
same marginal distribution and, to begin with, satisfy the following 
assumptions:   
\begin{itemize}
\item[{\rm (i)}]\  ${\rm E} Y_0=\mu>0$;
\item[{\rm (ii)}]\  ${\rm E}(Y_0^{2})<\infty$.
\end{itemize}

In terms of the generic sequence $\{Y_j,\, j=0,1,2,\ldots\}$, with 
$t\geq 0$, we define
\begin{eqnarray}
S(t)&:=&\sum_{i=1}^{[t]}Y_i,\label{st}\\
N(t)&:=&\inf\{s\geq 1:\, S(s)>t\},\label{nt}\\
Q(t)&:=& S(t)+\mu N(\mu t)-2\mu t,\label{qt}
\end{eqnarray}
whose respective appropriately normalized versions will be used in studying 
partial sums, their renewals, Bahadur-Kiefer type processes when the random 
variables in the sequence $Y_i,\, i=0,1,\ldots$ are weakly or strongly 
dependent.

The research area of what has become known as Bahadur-Kiefer processes was 
initiated by Bahadur \cite{Bahadur1966} (cf. also Kiefer \cite{Kiefer1967})
who established an almost sure representation of i.i.d. random variables 
based sample quantiles in terms of their empiricals. Kiefer 
\cite{Kiefer1970} substantiated this work via studying the deviations between 
the sample quantile and its empirical processes. These three seminal papers 
have since been followed by many related further investigations (cf., e.g.,
Cs\"org\H o and R\'ev\'esz \cite{CsR1978}, [12, Chapter 5], 
Shorack \cite{Shorack}, Cs\"org\H o \cite{Cs1983}, Deheuvels and Mason
\cite{DM}, \cite{DM2}, Deheuvels \cite{D921}, \cite{D922}, Cs\"org\H o and 
Horv\'ath [7, Chapters 3-6], Cs\"org\H o and Szyszkowicz \cite{CsSz}, 
and references in these works). 

It follows from the results of Kiefer \cite{Kiefer1970}, and also from 
Vervaat \cite{Verv1}, \cite{Verv2} as spelled out in Cs\'aki {\it et al.}
\cite{CsCsRS}, that the original i.i.d. based Bahadur-Kiefer process
cannot converge weakly to any non-degenerate random element of the 
$D[0,1]$ function space. On the other hand, Cs\"org\H o {\it et al.} 
\cite{CsSzW} showed the opposite conclusion to hold true for long-range 
dependence based Bahadur-Kiefer processes. For an illustration and discussion 
of this conclusion, we refer to the Introduction and Corollary 1.2 of Cs\'aki
{\it et al.} \cite{CsCsK}. For further results along these lines, we refer to 
Cs\"org\H o and Kulik \cite{CsK1}, \cite{CsK2}.

The study of the almost sure asymptotic behaviour of Bahadur-Kiefer type 
processes for sums and their renewals in the i.i.d. case was initiated by 
Horv\'ath \cite{Horv}, Deheuvels and Mason \cite{DM}, and augmented by
further references and results as in Cs\"org\H o and Horv\'ath [7,
Chapter 2]. 

Vervaat \cite{Verv1}, \cite{Verv2} initiated the study of limit theorems in 
general for processes with a positive drift and their inverses. For results 
on the asymptotic behaviour of integrals of Bahadur-Kiefer type processes 
for sums and their renewals, the so-called Vervaat processes, we refer to 
\cite{CsCsRS} in the i.i.d. case, \cite{CsakiCsorgoKulik} in the weakly 
dependent case, and \cite{CsCsK} in the strongly dependent case.

Back to the topics of this paper on Bahadur-Kiefer type processes for sums 
and their renewals, the forthcoming Section 2 is concerned with 
the weakly dependent case, and Section 3 concludes results in terms of 
long-range dependent sequences of random variables. Both of these sections
contain the relevant proofs as well.   

\section{Weakly dependent case}
 \renewcommand{\thesection}{\arabic{section}} \setcounter{equation}{0}

In this section we deal with weakly dependent random variables based 
Bahadur-Kiefer type processes. First we summarize the 
main results in the case when $Y_i$ are i.i.d. random variables with finite 
4-th moment. 

\medskip\noindent
{\bf Theorem A} {\it Assume that $\{Y_i,\, i=0,1,\ldots\}$ are i.i.d. random 
variables with ${\rm E} Y_0=\mu>0$, ${\rm E} (Y_0-\mu)^2=\sigma^2>0$, and 
${\rm E} Y_0^4<\infty$. Then we have
\begin{equation}
Q(T)=\sigma
\left(W(T)-W\left(T-\frac{\sigma}{\mu}W(T)\right)\right)+o_{\AS}(T^{1/4}),
\quad {\rm as}\ \  T\to\infty,
\label{qtinv}
\end{equation}
\begin{equation}
\limsup_{T\to\infty}\frac{\sup_{0\leq t\leq T}|Q(t/\mu)|}
{(T\log\log T)^{1/4}(\log T)^{1/2}}=\frac{2^{1/4}\sigma^{3/2}}{\mu^{3/4}},
\quad {\rm a.s.},
\label{qtlimsup}
\end{equation}
\begin{equation}
\lim_{T\to\infty}\frac{\sup_{0\leq t\leq T}|Q(t/\mu)|}
{(\log T)^{1/2}(\sup_{0\leq t\leq T}|\mu N(t)-t|)^{1/2}}=
\frac{\sigma}{\mu^{1/2}},
\quad {\rm a.s.},
\label{qtlim}
\end{equation}
\begin{equation}
\lim_{T\to\infty}{\rm P}(T^{-1/4}|Q(T/\mu)|\leq y)=
2\int_{-\infty}^\infty \Phi(y\mu^{3/4}\sigma^{-3/2}|x|^{-1/2})
\varphi(x)\, dx-1,
\label{qtl}
\end{equation}
where $\Phi$ is the standard normal distribution function and $\varphi$ is 
its density.
}

\medskip
We note that (\ref{qtinv}) and (\ref{qtl}) are due to Cs\"org\H o and 
Horv\'ath \cite{CsH}, (\ref{qtlimsup}) is due to Horv\'ath \cite{Horv} and 
(\ref{qtlim}) is due to Deheuvels and Mason \cite{DM}. All these results 
can be found in \cite{CsH}.

For the case of i.i.d. random variables when the 4-th moment does not exist, 
we refer to Deheuvels and Steinebach \cite{DS}.  

In this section we assume that $S(t)$ can be approximated by a standard 
Wiener process as follows.

\noindent
{\bf Assumption A} {\it On the same probability space there exist a sequence
$\{Y_i,\ i=0,1,2,\ldots\}$ of random variables, with the same marginal 
distribution, satisfying assumptions {\rm (i)} and {\rm (ii)}, and a 
standard Wiener process $W(t),\ t\geq 0$, such that
\begin{equation}
\sup_{0\leq t\leq T}|S(t)-\mu t-\sigma W(t)|=O_{\AS}(T^\beta)
\label{appr}
\end{equation}
almost surely, as $T\to\infty$, with $\sigma>0$, 
where $S(t)$ is defined by {\rm (\ref{st})} and $\beta<1/4$.}

In the case of $1/4\leq \beta<1/2$, there is a huge literature on strong 
approximation of the form (\ref{appr}) for weakly dependent random variables 
$\{Y_i\}$. The case $\beta<1/4$ is treated in Berkes {\it et al.} 
\cite{BerkesLiuWu}, where Koml\'os-Major-Tusn\'ady 
\cite{KomlosMajorTusnady1975} type strong approximations as in 
(\ref{appr}) are proved under fairly general assumptions of dependence.  
For exact statements of, and conditions for, strong approximations that
yield (\ref{appr}) to hold true for the partial sums as in Assumption A,
we refer to \cite{BerkesLiuWu}.

\begin{thm}
Under {\rm Assumption A} all the results {\rm (\ref{qtinv}), (\ref{qtlimsup}),
(\ref{qtlim}) and (\ref{qtl})} in {\rm Theorem A} remain true.
\end{thm}

\noindent
{\bf Proof}. In fact, we only have to prove (\ref{qtinv}), for the other 
results follow from the latter. It follows from \cite{CsH}, Theorem 1.3 
on p. 37, that under Assumption A we have
$$
\limsup_{T\to\infty}\frac{\sup_{0\leq t\leq T}
\left|\frac{t}{\mu}-N(t)-\frac{\sigma}{\mu}W(t/\mu)\right|}
{(T\log\log T)^{1/4}(\log T)^{1/2}}
=2^{1/4}\sigma^{3/2}\mu^{-7/4}\quad {\rm a.s.}
$$
and also
$$
\sup_{0\leq t\leq T}|\mu t-\mu S(N(\mu t))|=O_{\AS}(T^\beta)
$$
as $T\to\infty$. Hence, as $T\to\infty$, we arrive at
$$
Q(T)=S(T)+\mu N(\mu T)-2\mu T=S(T)-\mu T-\left(S(N(\mu T))-\mu N(\mu T)\right)
+O_{\AS}(T^\beta)
$$
$$
=\sigma(W(t)-W(N(\mu T)))+O_{\AS}(T^\beta)=
$$
$$
\sigma
\left(W(T)-W\left(T-\frac{\sigma}{\mu}W(T)\right)\right)+o_{\AS}(T^{1/4}),
$$
i.e., having (\ref{qtinv}) as desired.
\qed

\section{Strongly dependent case}
\renewcommand{\thesection}{\arabic{section}} \setcounter{equation}{0}
 \setcounter{thm}{0} \setcounter{lemma}{0}

In this section we deal with long range (strongly) dependent sequences, 
based on moving averages as defined by
\begin{equation}
\eta_j=\sum_{k=0}^\infty \psi_k \xi_{j-k},\quad j=0,1,2,\ldots,
\label{linear}
\end{equation}
where $\{\xi_k,\, -\infty<k<\infty\}$ is a double sequence of independent
standard normal random variables, and the sequence of weights
$\{\psi_k, \, k=0,1,2,\ldots\}$ is square summable. Then
${\rm E}(\eta_0)=0$, ${\rm E}(\eta_0^2)=\sum_{k=0}^\infty\psi_k^2=:
\sigma^2$ and, on putting $\tilde\eta_j=\eta_j/\sigma$, $\{\tilde\eta_j,\,
j=0,1,2,\ldots\}$ is a stationary Gaussian sequence with
${\rm E}(\tilde\eta_0)=0$ and ${\rm E}(\tilde\eta_0^2)=1$. If
$\psi_k\sim k^{-(1+\alpha)/2}\ell(k)$ with a slowly varying function,
$\ell(k)$, at infinity, then 
${\rm E}(\eta_j\eta_{j+n})\sim b_\alpha n^{-\alpha}\ell^2(n)$,
where the constant $b_\alpha$ is defined by
$$
b_\alpha=\int_0^\infty x^{-(1+\alpha)/2}(1+x)^{-(1+\alpha)/2}\, dx.
$$
Now let $G(\cdot)$ be a real valued Borel measurable function,
and define the subordinated sequence $Y_j=G(\tilde\eta_j)$, $j=0,1,2,\ldots$.
We assume throughout that $J_1:={\rm E}(G(\tilde\eta_0)\tilde\eta_0)\not=0$.
We say in this case that the Hermite rank of the function $G(\cdot)$ is 
equal to 1 (cf. Introduction of \cite{CsCsK}).

For $1/2<H<1$ let $\{W_H(t),\, t\geq 0\}$  be a fractional Brownian
motion (fbm), i.e., a mean-zero Gaussian process with covariance
\begin{equation}
{\rm E} W_H(s)W_H(t)=\frac{1}{2}(s^{2H}+t^{2H}-|s-t|^{2H}).
\label{covar}
\end{equation}

Based on a strong approximation result of Wang {\it et al.} 
\cite{WangLinGulati2003}, what follows next, was proved in Section 2 of
Cs\'aki {\it et al.} \cite{CsCsK}.

\medskip\noindent 
{\bf Theorem B}
{\it Let $\eta_j$ be defined by {\rm (\ref{linear})} with $\psi_k\sim
k^{-(1+\alpha)/2}$, $0<\alpha<1$, and put $\tilde\eta_j=\eta_j/\sigma$
with $\sigma^2:={\rm E}(\eta_0^2)=\sum_{k=0}^\infty \psi_k^2$. Let $G(\cdot)$
be a function whose Hermite rank is $1$, and put $Y_j=G(\tilde\eta_j)$,
$j=0,1,2,\ldots$. Furthermore, let $\{S(t),\, t\geq 0\}$
be as in {\rm (\ref{st})} and assume condition {\rm (ii)}. Then, on
an appropriate probability space for the sequence
$\{Y_j=G(\tilde\eta_j),\, j=0,1,\ldots\}$, one can construct a fractional
Brownian motion $W_{1-\alpha/2}(\cdot)$ such that, as $T\to\infty$, we have
\begin{equation}\label{eq:3/1}
\sup_{0\le t\le T}\left|S(t)-\mu t-
\frac{J_1 \kappa_{\alpha}}{\sigma}W_{1-\alpha/2}(t)\right|=
o_{\AS}(T^{\gamma/2+\delta}),
\end{equation}
where $\mu={\rm E}(Y_0)$, 
\begin{equation}
\kappa_{\alpha}^2=
2\frac{\int_0^{\infty}x^{-(\alpha+1)/2}(1+x)^{-(\alpha+1)/2}\,
dx}{(1-\alpha)(2-\alpha)},
\label{kalpha}
\end{equation}
$\gamma=2-2\alpha$ for $\alpha<1/2$,
$\gamma=1$ for $\alpha\geq 1/2$ and $\delta>0$ is arbitrary}.

{\it Moreover, if we also assume condition {\rm (i)}, then, as $T\to\infty$, }
\begin{equation}
\sup_{0\leq t\leq T}\left |\mu N(\mu t)-\mu t+\frac{J_1
\kappa_\alpha}{\sigma}
W_{1-\alpha/2}(t)\right
|=o_{\AS}(T^{\gamma/2+\delta}+T^{(1-\alpha/2)^2+\delta}),
\label{eq:main2}
\end{equation}
{\it with $\gamma$ as right above, and arbitrary $\delta>0$.}

\medskip
Now, for use in the sequel, we state iterated logarithm results for 
fractional Brownian motion and its increments, which follows from 
Ortega's extension in \cite{Ortega1984} of Cs\"org\H o and R\'ev\'esz 
\cite{CsR1979}, [12, Section 1.2].

\medskip\noindent
{\bf Theorem C} {\it
For $T>0$ let $a_T$ be a nondecreasing function of $T$ such that
$0<a_T\leq T$ and $a_T/T$ is nonincreasing. Then
\begin{equation}
\limsup_{T\to\infty}\frac{\sup_{0\leq t\leq T-a_T}\sup_{0\leq s\leq a_T}
|W_{1-\alpha/2}(t+s)-W_{1-\alpha/2}(t)|}
{a_T^{1-\alpha/2}(2(\log T/a_T+\log\log T))^{1/2}}=1\quad {\rm a.s.}
\label{Ortega}
\end{equation}
If $\lim_{T\to\infty}(\log(T/a_T))/(\log\log T)=\infty$, then we have $\lim$
instead of $\limsup$ in {\rm (\ref{Ortega})}.} 

\medskip
First we give an invariance principle for $Q(T)$ defined by (1.3) if
$\gamma/2<(1-\alpha/2)^2$, which corresponds to the i.i.d. case when the 
forth moment exists. Equivalently, we assume that
\begin{equation}
0<\alpha<2-\sqrt{2}.
\label{alfa}
\end{equation}
Note that in (\ref{qtinv2}) below, the random time argument of 
$W_{1-\alpha/2}$ is strictly positive for large enough $T$ with probability 1. 
So, without loss of generality, we may define 
$W_{1-\alpha/2}(T-u)=0$ if $u>T$.
 
\begin{thm} Under the conditions of {\rm Theorem B}, including {\rm (i)}
and {\rm (ii)}, assuming {\rm (\ref{alfa})}, as $T\to\infty$, we have
$$
Q(T)=
\frac{J_1\kappa_\alpha}{\sigma}(W_{1-\alpha/2}(T)-W_{1-\alpha/2}(N(\mu T))
+o_{\AS}(T^{\gamma/2+\delta})
$$
\begin{equation}
=\frac{J_1\kappa_\alpha}{\sigma}
\left(W_{1-\alpha/2}(T)-W_{1-\alpha/2}\left(T-
\frac{J_1\kappa_\alpha}{\sigma\mu}W_{1-\alpha/2}(T)\right)\right)
+o_{\AS}(T^{\gamma/2+\delta}).
\label{qtinv2}
\end{equation}
\end{thm}

\noindent {\bf Proof.} Put $c=J_1\kappa_\alpha/\sigma$. Then
$$
Q(T)=S(T)-\mu T+\mu N(\mu T)-\mu T
$$
$$
=cW_{1-\alpha/2}(T)+
o_{\AS}(T^{\gamma/2+\delta})+\mu (N(\mu T)-T).
$$
But
$$
\mu(T-N(\mu T))=S(N(\mu T))-\mu N(\mu T)+\mu T-S(N(\mu T))
$$
$$
=cW_{1-\alpha/2}(N(\mu T))+o_{\AS}((N(\mu T))^{\gamma/2+\delta})+
\mu T-S(N(\mu T)),
$$
and using (\ref{eq:main2}) and Theorem C, we have
$$
cW_{1-\alpha/2}(N(\mu T))=cW_{1-\alpha/2}\left(T-
\frac{c}{\mu}W_{1-\alpha/2}(T)+o_{\AS}(T^{\gamma/2+\delta}
+T^{(1-\alpha/2)^2+\delta})\right)
$$
$$
=cW_{1-\alpha/2}\left(T-\frac{c}{\mu}W_{1-\alpha/2}(T)\right)+
o_{\AS}(T^{(\gamma/2+\delta)(1-\alpha/2)}+T^{(1-\alpha/2)^3}).
$$
On the other hand (cf. \cite{CsCsK}), $N(\mu T)=O_{\AS}(T)$ and
$$
\mu T-S(N(\mu T))=o_{\AS}(T^{\gamma/2+\delta}).
$$
Since $(1-\alpha/2)^3\leq \gamma/2< (1-\alpha/2)^2$, this 
dominates all the other remainder terms
in the proof. Thus the proof of Theorem 3.1 is now complete.
\qed

The proof of Theorem 3.1 also yields the following result.
\begin{proposition}
As $T\to\infty$,
$$
\mu T-\mu N(\mu T))=\frac{J_1\kappa_\alpha}{\sigma}
W_{1-\alpha/2}\left(T-\frac{J_1\kappa_\alpha}{\sigma\mu}W_{1-\alpha/2}(T)
\right)+o_{\AS}(T^{\gamma/2+\delta}).
$$
\end{proposition}

Now we are to give a limsup result for $Q(\cdot)$. For this we need a 
Strassen-type functional law of the iterated logarithm for fbm, due to 
Goodman and Kuelbs \cite{GK1991}.

\newpage
\medskip\noindent{\bf Theorem D} {\it Let
$$
{\bf K}=\{T_H g(t),\, 0\leq t\leq 1,\, 
\int_{-\infty}^1 g^2(u)\, du\leq 1\},
$$ 
where
$$
T_H g(t)=\frac{1}{k_H}\int_0^t (t-u)^{H-1/2}g(u)\, du
+\frac{1}{k_H}\int_{-\infty}^0(t-u)^{H-1/2}-(-u)^{H-1/2})g(u)\, du,
$$
and
$$
k_H^2=\int_{-\infty}^0((1-s)^{H-1/2}-(-s)^{H-1/2})^2\, ds
+\int_0^1(1-s)^{2H-1}\, ds.
$$
Then, almost surely, ${\bf K}$ is the set of limit points of the net of
stochastic processes
\begin{equation}
\frac{W_H(nt)}{(2n^{2H}\log\log n)^{1/2}},\, 0\leq t\leq 1,
\label{flil}
\end{equation}
as $n\to\infty$.}

\begin{thm} Under the conditions of {\rm Theorem 3.1}, we have
\begin{equation}
\limsup_{T\to\infty}
\frac{|Q(T)|}{T^{(1-\alpha/2)^2}(\log\log T)^{1/2-\alpha/4}(\log T)^{1/2}}
=\frac{2^{1-\alpha/4} (J_1\kappa_\alpha)^{2-\alpha/2}}
{\sigma^{2-\alpha/2}\mu^{1-\alpha/2}}\qquad {\rm a.s.}
\label{qtlimsup2}
\end{equation}
\end{thm}

\noindent
{\bf Proof}. It follows from Theorem C that
$$
|W_{1-\alpha/2}(T)|\leq (1+\delta)T^{1-\alpha/2}(2\log\log T)^{1/2}
$$
with probability 1 for any $\delta>0$ if $T$ is large enough. Hence, applying 
Theorem C with $a_T=(1+\delta)c/\mu T^{1-\alpha/2}(2\log\log T)^{1/2}$,  
$c=J_1\kappa_\alpha/\sigma$, we obtain
$$
c\sup_{|s|\leq a_T}|W_{1-\alpha/2}(T)-W_{1-\alpha/2}(T-s)|\leq
c(1+\delta)a_T^{1-\alpha/2}(2\log T)^{1/2},
$$
almost surely for large enough T. Since $\delta>0$ is arbitrary, 
we obtain the upper bound in (\ref{qtlimsup2}).

To obtain the lower bound, we follow the proof in the i.i.d. case, given in 
Cs\"org\H o and Horv\'ath \cite{CsH}. On choosing 
\begin{displaymath}
g(s)=\left\{
\begin{array}{ll}
\frac{1}{k_H}((1-s)^{H-1/2}-(-s)^{H-1/2}),
&\quad s\leq 0,\\
\frac{1}{k_H}(1-s)^{H-1/2}, &\quad 0<s\leq 1,
\end{array}
\right.
\end{displaymath}
in Theorem D, we have
$$
f(t)=\frac{1}{k_H}\int_{-\infty}^0((t-s)^{H-1/2}-(-s)^{H-1/2})g(s)\, ds
+\frac{1}{k_H}\int_0^t(t-s)^{H-1/2}g(s)\, ds.
$$
It can be seen that $\int_{-\infty}^1 g^2(s)\, ds=1$, and 
$\{f(t),\, \, 0\leq t\leq 1\}$ is a continuous increasing function with 
$f(0)=0$, $f(1)=1$, and hence by Theorem D it is in ${\bf K}$. 
For $0<\delta<1$, on considering the function
\begin{displaymath}
g_\delta(s)=\left\{
\begin{array}{ll}
g(s), &\quad 0\leq s\leq 1-\delta,\\ 
0, &\quad 1-\delta\leq s\leq 1,
\end{array}
\right.
\end{displaymath}
we define
\begin{displaymath}
f_\delta(t)=\left\{
\begin{array}{ll}
f(t), &\quad 0\leq t\leq 1-\delta,\\
f(1-\delta), &\quad 1-\delta\leq t\leq 1.
\end{array}
\right.
\end{displaymath}
Then it can be seen that the latter function is in ${\bf K}$, and hence it 
is a limit function of the net of stochastic processes as in (\ref{flil}). 
It follows that there is a sequence $T_k$ of random variables such that,
in our context,
$$
\lim_{k\to\infty}\sup_{0\leq t\leq 1}
\left|\frac{W_{1-\alpha/2}(T_kt)}
{T_k^{1-\alpha/2}(2\log\log T_k)^{1/2}}-f_\delta(t)\right|=0.
$$
Using Theorem C with $a_T=f(1-\delta)c/\mu T^{1-\alpha/2}(2\log\log T)^{1/2}$,
we get
$$
\lim_{T\to\infty}\frac{\sup_{T(1-\delta)\leq t\leq T}
c|W_{1-\alpha/2}(t+a_T)-W(t)|}{ca_T^{1-\alpha/2}(2\log T)^{1/2}}=1
\qquad {\rm a.s.}
$$

Since $\delta$ is arbitrary, and 
$\lim_{\delta\to 0}f_\delta(t)=f(t)$, the lower bound  
follows as in Cs\"org\H o and Horv\'ath \cite{CsH}, p. 28. This completes the 
proof of Theorem 3.2.
\qed

\medskip
Next we give the limiting distribution of $Q(T)$.
\begin{thm} Under the conditions of {\rm Theorem 3.1}, we have
\begin{equation}
\lim_{T\to\infty}{\rm P}\left(Q(T)T^{-(1-\alpha/2)^2}\leq y\right)= 
\int_{-\infty}^{\infty}\varphi(x)
\Phi\left(\frac{y \sigma^{2-\alpha/2}\mu^{1-\alpha/2}}
{|x|^{1-\alpha/2}(J_1\kappa_\alpha)^{2-\alpha/2}}\right)\, dx.
\label{limit}
\end{equation}
\end{thm}

\noindent{\bf Proof.}
According to Theorem 3.1 we have to determine the limiting distribution of
$$
c\left(W_{1-\alpha/2}(T)-W_{1-\alpha/2}\left(T-\frac{c}{\mu}
W_{1-\alpha/2}(T)\right)\right),
$$
where $c=J_1\kappa_\alpha/\sigma$. Via the scaling property of fbm, i.e.,
$$
\widetilde W(v):=T^{-1+\alpha/2}W_{1-\alpha/2}(Tv), \quad v\geq 0,
$$
is also an fbm with parameter $1-\alpha/2$. So we have to determine the 
limiting distribution of 
$$c\left(\widetilde W(1)-\widetilde 
W(1-c_1T^{-\alpha/2}\widetilde W(1))\right),$$ 
as $T\to\infty$, where $c_1=J_1\kappa_\alpha/(\sigma\mu)$.

For $u>0$, the joint distribution of $\widetilde W(1)$, $\widetilde W(u)$
is bivariate normal with density
$$
\frac{1}{2\pi\sigma_1\sigma_2\sqrt{1-r^2}}
\exp\left(-\frac{1}{2(1-r^2)}
\left(\frac{x^2}{\sigma_1^2}-2r\frac{xy}{\sigma_1\sigma_2}+
\frac{y^2}{\sigma_2^2}\right)\right),
$$
where $\sigma_1^2={\rm E}(W_{1-\alpha/2}^2(1))=1$,
$\sigma_2^2={\rm E}(W_{1-\alpha/2}^2(u))=u^{2-\alpha}$ and
$$
r=\frac{1+u^{2-\alpha}-|1-u|^{2-\alpha}}{2\sigma_1\sigma_2}.
$$
Now consider the conditional density 
$$
{\rm P}(\widetilde W(u)\in dz|\widetilde W(1)=x)=
\frac{1}{\sigma_2\sqrt{1-r^2}}\varphi\left(\frac{z-r\sigma_2 x}
{\sigma_2\sqrt{1-r^2}}\right)\, dz,
$$
where $u=1-c_1xT^{-\alpha/2}$.

So the density function of $\widetilde W(1)-\widetilde W(u)$ is equal to
$$
{\rm P}(\widetilde W(1)-\widetilde W(u)\in dY)=
\int_{-\infty}^{T^{\alpha/2}/c_1}\frac{1}{\sigma_2\sqrt{1-r^2}}\varphi(x)
\varphi\left(\frac{x-Y-r\sigma_2x}{\sigma_2\sqrt{1-r^2}}\right)\, dx\, dY
$$
and hence its distribution function is 
$$
{\rm P}(\widetilde W(1)-\widetilde W(u)\leq Z)=
\int_{-\infty}^{T^{\alpha/2}/c_1}\varphi(x)
\Phi\left(\frac{Z-x+r\sigma_2x}{\sigma_2\sqrt{1-r^2}}\right)\, dx,
\quad -\infty<Z<\infty.
$$ 
It can be seen that, as $T\to\infty$,
$$
\sigma_2\sqrt{1-r^2}\sim \frac{|c_1x|^{1-\alpha/2}}{T^{\alpha/2-\alpha^2/4}},
$$
$$
\frac{x-xr\sigma_2}{\sigma_2\sqrt{1-r^2}}=O(T^{-\alpha/2+\alpha^2/4}).
$$
Hence, as $T\to\infty$,
$$
{\rm P}(\widetilde W(1)-\widetilde W(u)\leq Z)\sim
\int_{-\infty}^{T^{\alpha/2}/c_1}\varphi(x)
\Phi\left(\frac{ZT^{\alpha/2-\alpha^2/4}}{|c_1x|^{1-\alpha/2}}\right)\, dx.
$$
Putting $Z=yT^{\alpha^2/4-\alpha/2}/c$, and taking the limit $T\to\infty$, we 
finally obtain (\ref{limit}). 
\qed

\bigskip
\noindent
{\bf Acknowledgement} We wish to thank two referees for their careful 
reading of, and constructive remarks on, our manuscript. Research supported 
by an NSERC Canada Discovery Grant at Carleton University, Ottawa and by the 
Hungarian National Foundation for Scientific Research, No. K108615.


\begin{thebibliography}{9}

\bibitem{Bahadur1966} 
Bahadur, R.R.: A note on quantiles in large 
samples. \textit{Ann. Math. Statist.} \textbf{37}, 577--580 (1966)

\bibitem{BerkesLiuWu} 
Berkes, I., Liu, W.D. and Wu, W.B.: 
Koml\'os-Major-Tusn\'ady approximation under dependence. 
\textit{Ann. Probab.} \textbf{42}, 794--817 (2014)

\bibitem{CsakiCsorgoKulik} 
Cs{\'a}ki, E., Cs{\"o}rg{\H{o}}, M. and Kulik, R.:
On {V}ervaat processes for sums and renewals in weakly dependent
cases. In: \textit{{D}ependence in {P}robability, {A}nalysis and {N}umber
{T}heory}. A {V}olume in {M}emory of {W}alter {P}hilipp. Berkes {\it et al.},
ed., Kendrick Press, Heber City, UT., pp. 145--156 (2010)

\bibitem{CsCsK} 
Cs\'aki, E., Cs\"org\H o, M. and Kulik, R.:
Strong approximations for long memory sequences based partial sums,
counting and their Vervaat processes. Submitted. arXiv:math.PR1302.3740
(2013)

\bibitem{CsCsRS} 
Cs\'aki, E., Cs\"org\H o, M., Rychlik, Z. and Steinebach, J.:  
On Vervaat and Vervaat-error-type processes for partial sums
and renewals. \textit{J. Statist. Plann. Inf.} \textbf{137}, 953--966
(2007)

\bibitem{Cs1983}
Cs\"org\H o, M.: \textit{Quantile Processes with Statistical 
Application}. CBMS-NSF Regional Conference Series in Applied Mathematics
\textbf{42}, SIAM, Philadelphia (1983)

\bibitem{CsH}
Cs\"org\H o, M. and Horv\'ath, L.: \textit{Weighted Approximations
in Probability and Statistics}. Wiley, Chichester (1993)

\bibitem{CsK1}
Cs\"org\H o, M. and Kulik, R.:  
Reduction principles for quantile and Bahadur-Kiefer processes of long-range 
dependent sequences. \textit{Probab. Th. Rel. Fields} \textbf{142} 339--366
(2008)

\bibitem{CsK2}
Cs\"org\H o, M. and Kulik, R.:  
Weak convergesnce of Vervaat and Vervaat error processes of long-range 
dependent sequences. \textit{J. Theoret. Probab.} \textbf{21}, 672--686
(2008) 

\bibitem{CsR1978}
Cs{\"o}rg{\H{o}}, M. and R\'ev\'esz, P.:  
Strong Approximations of the quantile process. \textit{Ann. Statist.} 
\textbf{6}, 882--894 (1978)

\bibitem{CsR1979}
Cs{\"o}rg{\H{o}}, M. and R\'ev\'esz, P.:  
How big are the increments of a Wiener process? \textit{Ann. Probab.} 
\textbf{7}, 731--737 (1979) 

\bibitem{CsR1981}
Cs{\"o}rg{\H{o}}, M. and R\'ev\'esz, P.:
\textit{Strong Approximations in Probability and Statistics.}
Academic Press, New York (1981)

\bibitem{CsSz} 
Cs\"org\H o, M., Szyszkovicz, B.: 
Sequential quantile and Bahadur-Kiefer processes. \textit{Order Statistics: 
Theory and Methods}. Handbook of Statistics, \textbf{16},  
North-Holland, Amsterdam, pp. 631--688 (1998)

\bibitem{CsSzW}
Cs\"org\H o, M., Szyszkowicz, B. and Wang, L.H.: 
Strong invariance principles for sequential Bahadur-Kiefer and Vervaat error 
processes of long-range dependent sequences. \textit{Ann. Statist.} 
\textbf{34}, 1013--1044 (2006).
Correction: \textit{Ann. Statist.} \textbf{35}, 2815--2817 (2007) 

\bibitem{D921}
Deheuvels, P.: 
Pointwise Bahadur-Kiefer-type theorems I. 
\textit{Probability Theory and Applications.} 235-255, Math. Appl., 80, Kluwer 
Acad. Publ., Dordrecht. pp. 235--255 (1992)

\bibitem{D922}
Deheuvels, P.: 
Pointwise Bahadur-Kiefer-type theorems II. 
\textit{Nonparametric Statistics and Related Topics (Ottawa, ON)},
North-Holland, Amsterdam. pp. 331--345 (1992)

\bibitem{DM}
Deheuvels, P. and Mason, D.M.:  
Bahadur-Kiefer-type processes. \textit{Ann. Probab.} \textbf{18}, 669--697
(1990)

\bibitem{DM2}
Deheuvels, P. and Mason, D.M.:  
A functional LIL approach to pointwise Bahadur-Kiefer theorems. 
\textit{Probability in Banach Spaces,} 
\textbf{8} (Brunswick, ME, 1991) Progr. Probab., 30, Birkh\"auser 
Boston, Boston, MA. pp. 255--266 (1992)

\bibitem{DS}
Deheuvels, P. and Steinebach, J.: 
On the limiting behavior of the Bahadur-Kiefer statistic for partial sums 
and renewal processes when the fourth moment does not exist. 
\textit{Statist. Probab. Lett.} \textbf{13}, 179--188 (1992)

\bibitem{GK1991}
Goodman, V. and Kuelbs, J.: 
Rates of clustering for some Gaussian self-similar processes. 
\textit{Probab. Th. Rel. Fields} \textbf{88}, 47--75 (1991)

\bibitem{Horv}
Horv\'ath, L.: 
Strong approximations of renewal processes.
\textit{Stoch.\ Process.\ Appl.} \textbf{18} 127--138 (1984)

\bibitem{Kiefer1967}
Kiefer, J.: 
On {B}ahadur's representation of sample quantiles.
\textit{Ann. Math. Statist.} \textbf{38}, 1323--1342 (1967)

\bibitem{Kiefer1970}
Kiefer, J.:  
Deviations between the sample quantile process and the 
sample df. \textit{Nonparametric Techniques in Statistical Inference},
Cambridge Univ. Press, London, pp. 299--319 (1970) 

\bibitem{KomlosMajorTusnady1975}
Koml{\'o}s, J., Major, P. and Tusn{\'a}dy, G.: 
An approximation of partial sums of independent {${\rm RV}$}'s and
the sample {${\rm DF}$}. {I}. \textit{Z. Wahrsch. Verw. Gebiete} 
\textbf{32}, 111--131 (1975)

\bibitem{Ortega1984}
Ortega, J.: 
On the size of the increments of nonstationary 
{G}aussian processes. \textit{Stoch. Process. Appl.} 
\textbf{18}, 47--56 (1984)

\bibitem{Shorack}
Shorack, G.R.: 
Kiefer's theorem via the Hungarian construction.
\textit{Z. Wahrsch. Verw. Gebiete} \textbf{61} 369--373 (1982)

\bibitem{Verv1}
Vervaat, W.: 
\textit{Success Epochs in Bernoulli Trials: with
Applications to Number Theory.} Mathematical Centre Tracts \textbf{42}
(second edition in 1977). Matematisch Centrum, Amsterdam (1972)

\bibitem{Verv2}
Vervaat, W.: 
Functional central limit theorems for processes with positive drift and 
their inverses. \textit{Z.\ Wahrsch.\ Verw.\ Gebiete} 
\textbf{23}, 245--253 (1972)

\bibitem{WangLinGulati2003}
Wang, Q., Lin, Y-X. and Gulati, C.M.: 
Strong approximation for long memory processes with applications.
\textit{J. Theoret. Probab.} \textbf{16} 377--389 (2003)

\end{thebibliography}
\end{document}